\documentclass[default,a4]{sn-jnl}% Default
\jyear{2021}%

\usepackage{lmodern}   % 数学公式字体

\usepackage{autobreak}       %e autobreak package

 \usepackage[capitalise]{cleveref}%\usepackage[nameinlink]{cleveref} %  宏包的\cref 类似于\cite引用多个文献
 \usepackage[shortlabels]{enumitem}  % 设定排列环境的 格式, 可以缩写简略形式
  \usepackage{extarrows}  % 定义箭头  \xlongequal  \xLongleftrightarrow
\numberwithin{equation}{section} %公式按章节编号
 \allowdisplaybreaks[4]   % 使数学公式自动分页
\newtheorem{definition}{Definition} [section]            % 斜体
\newtheorem{theorem}{Theorem}[section]            % 斜体
\newtheorem{lemma}{Lemma} [section]%单独编号,与节有关
\newtheorem{proposition}{Proposition}[section]
\newtheorem{remark}{Remark}
\newenvironment{refproof}[1][]     % 若少了[],后边就需要改为{XXX} 定义 Proof of Theorem X
  {\renewcommand{\proofname}{\bf #1}%
   \begin{proof}}
  {\end{proof}}
\def\loc{{\mathrm{loc}}}

\def\dchi{\scalebox{1.2}{$\chi$}}

\def\dint{\displaystyle\int}

\DeclareMathOperator*{\bmo}{BMO}

\DeclareMathOperator*{\essinf}{ess\, inf}
\DeclareMathOperator*{\esssup}{ess\, sup}
\DeclareMathOperator*{\lip}{Lip}

\newcommand{\mathd}{\mathrm{d}}

\begin{document} %\bfseries

\title[Commutators of  Maximal Function on Stratified Lie Groups]{{ Characterization of Lipschitz Functions  via the Commutators of  Maximal Function on Stratified Lie Groups }}


\author[1]{\small \fnm{Jianglong} \sur{Wu}}%\email{jl-wu@163.com}

\author*[2]{\small\fnm{Wenjiao} \sur{Zhao}}%\email{wenjiaozhao@163.com}


\affil[1]{ \footnotesize \orgdiv{Department of Mathematics}, \orgname{Mudanjiang Normal University}, \orgaddress{ \city{Mudanjiang}, \postcode{157011},\country{China}}}

%\affil*[2]
\affil*[2]{ \footnotesize \orgdiv{School of Mathematics}, \orgname{Harbin Institute of Technology}, \orgaddress{ \city{Harbin}, \postcode{150001},   \country{China}}}

%\affil[3]{\orgdiv{Department}, \orgname{Organization}, \orgaddress{\street{Street}, \city{City}, \postcode{610101}, \state{State}, \country{Country}}}


\abstract{
In this paper, the main aim  is to consider the boundedness of  the   Hardy-Littlewood maximal commutator $M_{b}$ and the nonlinear commutator   $[b, M]$  on the Lebesgue spaces and Morrey spaces over some stratified Lie group  $\mathbb{G}$ when $b$ belongs to the Lipschitz space, by which some new characterizations of the Lipschitz spaces on Lie group  are given.
}




\maketitle
%==========================
\footnote{*Corresponding author \\
Email addresses: wenjiaozhao@163.com}
%==========================
\section{Introduction and main results}
\label{sec:introduction}

Stratified groups appear in quantum physics and many parts of mathematics, including  several complex variables, Fourier analysis, geometry, and topology \cite{folland1982hardy,varopoulos2008analysis}.
The geometry structure of stratified groups is so good that it inherits a lot of analysis properties from the Euclidean spaces \cite{stein1993harmonic,grafakos2009modern}.
Apart from this, the difference between the geometry structures of Euclidean spaces and stratified groups makes the study of function spaces on them more complicated.
However, many harmonic analysis problems  on stratified Lie groups deserve a further investigation since  most results   of the theory of Fourier transforms and distributions in Euclidean spaces  cannot yet be duplicated.

Let $T$ be the classical singular integral operator. The commutator $[b, T]$ generated by $T$ and a
suitable function $b$ is defined by
\begin{align} \label{equ:commutator-1}
 [b,T]f      & = bT(f)-T(bf).
\end{align}
It is  known that the commutators are intimately related to the regularity properties of the solutions of certain partial differential equations (PDE), see   \cite{difazio1993interior,bramanti1995commutators,rios2003lp}.

The first result for the commutator $[b,T]$ was established by  % Coifman, Rochberg and Weiss in
\citet{coifman1976factorization}, and the authors proved that $b\in \bmo(\mathbb{R}^{n})$ (bounded mean oscillation functions) if and only if the commutator \labelcref{equ:commutator-1} is  bounded on $L^{p}(\mathbb{R}^{n})$  for $1<p<\infty$.
 In 1978,  Janson \cite{janson1978mean} generalized the results in \cite{coifman1976factorization} to functions belonging to a Lipschitz functional space and gave some characterizations of the Lipschitz space $\dot{\Lambda}_{\beta}(\mathbb{R}^{n})$ via commutator \labelcref{equ:commutator-1}, and the author proved that $b\in \dot{\Lambda}_{\beta}(\mathbb{R}^{n})$ if and only if $[b,T]$ is bounded from $L^{p}(\mathbb{R}^{n})$ to $L^{q}(\mathbb{R}^{n})$ where $1<p<n/\beta$ and $1/p-1/q=\beta/n$ (see also  \cite{paluszynski1995characterization}).

 In  addition, using real interpolation techniques,  \citet{milman1990second} established a commutator result that applies to the Hardy-Littlewood maximal function as well as to a large class of nonlinear operators.
 In 2000, \citet{bastero2000commutators}    proved the necessary and sufficient conditions for the boundedness of  the  nonlinear commutator $[b,M]$ on $L^{p}$ spaces, and the similar problems for $[b,M_{\alpha}]$  were also studied by  \citet{zhang2009commutators}.
In 2017, \citet{zhang2017characterization} considered some new characterizations of the Lipschitz spaces via the boundedness of maximal commutator $M_{b}$ and the (nonlinear) commutator $[b, M]$ in Lebesgue spaces and Morrey spaces on Euclidean spaces.
In 2018, \citet{zhang2018commutators} gave necessary and sufficient conditions for the boundedness of the nonlinear commutator  $[b,M_{\alpha}]$  on Orlicz spaces when the symbol $b$ belongs to Lipschitz spaces, and obtained some new characterizations of non-negative Lipschitz functions.
And \citet{guliyev2022some} recently gave  necessary and sufficient conditions for the boundedness of  the fractional maximal commutators in the Orlicz spaces $L^{\Phi} (\mathbb{G})$ on  stratified Lie group  $\mathbb{G}$   when   $b$ belongs to $\bmo(\mathbb{G})$ spaces, and obtained some new characterizations for certain subclasses  of $\bmo(\mathbb{G})$ spaces.

Inspired by the above literature, the purpose of this paper is to   study the boundedness of the Hardy-Littlewood maximal commutator $M_{b}$ and the nonlinear commutator   $[b, M]$
in the Lebesgue spaces and Morrey spaces on some stratified Lie group  $\mathbb{G}$ when $b\in \dot{\Lambda}_{\beta}(\mathbb{G})$, by which some new characterizations of the Lipschitz spaces are given.

 Let  $f\in L_{\loc}^{1}(\mathbb{G})$,  the Hardy–Littlewood maximal function  $M$  is given by
%-----------------
\begin{align*}
%-------
M (f)(x) &= \sup_{B\ni x} \vert B \vert ^{-1}  \dint_{B} \vert f(y) \vert  \mathd y
%----------
\end{align*}
%------------
where the supremum is taken over all balls $B\subset \mathbb{G}$  containing $x$, and $\vert B \vert $  is the Haar measure of the $ \mathbb{G}$-ball $B$.
  And the maximal commutator $M_{b}$ generated by  the operator $M$  and a locally integrable function $b$ is defined by
%-----------------
\begin{align*}
%-------
M_{b} (f)(x) &= \sup_{B\ni x} \vert B \vert ^{-1}  \dint_{B} \vert b(x)-b(y) \vert  \vert f(y) \vert  \mathd y.
%----------
\end{align*}
%------------
On the other hand, similar to \labelcref{equ:commutator-1}, we can define the (nonlinear) commutator of the Hardy-Littlewood maximal function $M$  with a locally integrable function $b$ is defined by
%-----------------
\begin{align*}
%-------
 [b,M] (f)(x) &= b(x) M (f)(x) - M (bf)(x).
%----------
\end{align*}
%------------

Note that operators $M_{b}$ and $[b, M]$ essentially differ from each other. For example, $M_{b}$ is positive and sublinear, but $[b, M]$ is neither positive nor sublinear.

The first part of this paper is to study the boundedness of $M_{b}$  when the symbol $b$ belongs to a Lipschitz space. Some characterizations of the Lipschitz space via such commutator are given.

\begin{theorem} \label{thm:lipschitz-main-1}
Let $b$ be a locally integrable function and $0 <\beta <1$. Then the following statements are equivalent:
%-----------------------
\begin{enumerate}[label=(\arabic*),itemindent=2em]
 %-------------------------------------
  \item   $b\in \dot{\Lambda}_{\beta}(\mathbb{G})$.
 %-----------
    \label{enumerate:thm-lip-main-1-1}
%------------------
   \item $ M_{b} $ is bounded from $L^{p}(\mathbb{G})$ to $L^{q}(\mathbb{G})$ for all $p, q$ with $1<p<Q/\beta$ and $1/q = 1/p -\beta/Q$.
%-----------
    \label{enumerate:thm-lip-main-1-2}
%--------------------------------
\item $ M_{b} $ is bounded from $L^{p}(\mathbb{G})$ to $L^{q}(\mathbb{G})$ for some $p, q$ with $1<p<Q/\beta$ and $1/q = 1/p -\beta/Q$.
%-----------
    \label{enumerate:thm-lip-main-1-3}
%--------------------------------
   \item  $ M_{b} $   satisfies the weak-type $(1,Q/(Q-\beta))$ estimates, namely, there exists a positive constant $C$ such that
%------------------
\begin{align} \label{inequ:lip-main-1-4}
%-----------
   \big \lvert \{ x\in \mathbb{G}: M_{b}(f)(x) > \lambda \}  \big  \rvert \le C  \left( \lambda^{-1}    \|f  \|_{L^{1}(\mathbb{G})} \right)^{Q/(Q-\beta)}
%-----------------
\end{align}
%--------------
 holds for all   $\lambda > 0$.
%-----------
    \label{enumerate:thm-lip-main-1-4}
%-------------------
\end{enumerate}
%--------------------
\end{theorem}

%-----------------------
\begin{theorem} \label{thm:lipschitz-mb-main-2}
%-----------------------
Let  $b$ be a locally integrable function and $0 <\beta <1$. Suppose that $1 <p < Q/\beta$, $0<\lambda <Q-\beta p$.
%-----------------------
\begin{enumerate}[label=(\arabic*),itemindent=2em]
 %%-------------------------------------
  \item   If   $1/q=1/p-\beta/(Q-\lambda)$.  Then $b\in \dot{\Lambda}_{\beta}(\mathbb{G})$ if and only if $ M_{b} $  is bounded from  $L^{p,\lambda}(\mathbb{G})$ to $L^{q,\lambda}(\mathbb{G})$.
%-----------
    \label{enumerate:thm-lip-mb-main-2-1}
%--------------------------------
   \item If   $1/q=1/p-\beta/Q$ and $\lambda/p =\mu/q$.  Then $b\in \dot{\Lambda}_{\beta}(\mathbb{G})$ if and only if $ M_{b} $  is bounded from  $L^{p,\lambda}(\mathbb{G})$ to $L^{q,\mu}(\mathbb{G})$.
%-----------
    \label{enumerate:thm-lip-mb-main-2-2}
%--------------------------------
\end{enumerate}
%--------------------
\end{theorem}
%--------------------

The second part of this paper aims to study the mapping properties of the (nonlinear) commutator $[b, M]$ when  $b$ belongs to some Lipschitz space. To state our results, we recall the definition of the maximal operator with respect to a ball. For a fixed ball $B_{0}$, the Hardy-Littlewood maximal function with respect to $B_{0}$ of a function $f$ is given by
%-----------------
\begin{align*}
%-------
M_{B_{0}} (f)(x) &= \sup_{B_{0}\supseteq B\ni x} \vert B \vert ^{-1}  \dint_{B} \vert f(y) \vert  \mathd y,
%----------
\end{align*}
%------------
where the supremum is taken over all the balls $B$  with $B\subseteq B_{0}$ and $x\in B$.

%-----------------------
\begin{theorem} \label{thm:lipschitz-nonlinear-main-1}
%-----------------------
Let  $b$ be a locally integrable function and $0 <\beta <1$. Suppose that $1 <p < Q/\beta$ and $1/q=1/p-\beta/Q$.  Then the following statements are equivalent:
%-----------------------
\begin{enumerate}[label=(\arabic*),itemindent=2em]
 %%-------------------------------------
  \item   $b\in \dot{\Lambda}_{\beta}(\mathbb{G})$ and $b\ge 0$.
%-----------
    \label{enumerate:thm-lip-nonlinear-main-1-1}
%--------------------------------
   \item $[b,M ]$ is bounded from $L^{p}(\mathbb{G})$ to $L^{q}(\mathbb{G})$.
%-----------
    \label{enumerate:thm-lip-nonlinear-main-1-2}
%--------------------------------
   \item  There exists a constant $C > 0$ such that
%------------------
\begin{align}  \label{inequ:lip-nonlinear-main-1-3}
%-----------
  \sup_{B\ni x} \vert B \vert ^{-\beta/Q}  \left( \lvert B \rvert^{-1} \dint_{B}  \lvert b(x) -M_{B}(b)(x)  \rvert^{q} \mathd x \right)^{1/q}  \le C.
%-----------------
\end{align}
%-----------
    \label{enumerate:thm-lip-nonlinear-main-1-3}
%------------------------
\end{enumerate}
%--------------------
\end{theorem}
%--------------------

%-----------------------
\begin{theorem} \label{thm:lipschitz-nonlinear-main-3}
%-----------------------
Let  $b$ be a locally integrable function and $0 <\beta <1$. Suppose that $1 <p < Q/\beta$, $0<\lambda <Q-\beta p$ and  $1/q=1/p-\beta/(Q-\lambda)$.   Then the following statements are equivalent:
%-----------------------
\begin{enumerate}[label=(\arabic*),itemindent=2em]
 %%-------------------------------------
  \item   $b\in \dot{\Lambda}_{\beta}(\mathbb{G})$ and $b\ge 0$.
%-----------
    \label{enumerate:thm-lip-nonlinear-main-3-1}
%--------------------------------
   \item $[b,M ]$ is bounded from  $L^{p,\lambda}(\mathbb{G})$ to $L^{q,\lambda}(\mathbb{G})$.
%-----------
    \label{enumerate:thm-lip-nonlinear-main-3-2}
%------------------------
\end{enumerate}
%--------------------
\end{theorem}
%--------------------

%-----------------------
\begin{theorem} \label{thm:lipschitz-nonlinear-main-4}
%-----------------------
Let  $b$ be a locally integrable function and $0 <\beta <1$. Suppose that $1 <p < Q/\beta$, $0<\lambda <Q-\beta p$,  $1/q=1/p-\beta/Q$ and $\lambda/p =\mu/q$.   Then the following statements are equivalent:
%-----------------------
\begin{enumerate}[label=(\arabic*),itemindent=2em]
 %%-------------------------------------
  \item   $b\in \dot{\Lambda}_{\beta}(\mathbb{G})$ and $b\ge 0$.
%-----------
    \label{enumerate:thm-lip-nonlinear-main-4-1}
%--------------------------------
   \item $[b,M ]$ is bounded from  $L^{p,\lambda}(\mathbb{G})$ to $L^{q,\mu}(\mathbb{G})$.
%-----------
    \label{enumerate:thm-lip-nonlinear-main-4-2}
%------------------------
\end{enumerate}
%--------------------
\end{theorem}
%--------------------

This paper is organized as follows. In the next section, we recall some basic definitions and known results. In  \cref{sec:proof-mb}, we will prove  \cref{thm:lipschitz-main-1,thm:lipschitz-mb-main-2}.   \cref{sec:proof-nonlinear}  is devoted to proving \cref{thm:lipschitz-nonlinear-main-1,thm:lipschitz-nonlinear-main-3,thm:lipschitz-nonlinear-main-4}.

Throughout this paper, the letter $C$  always stands for a constant  independent of the main parameters involved and whose value may differ from line to line.

\section{Preliminaries and lemmas}
\label{sec:preliminary}

\subsection{Lie group $\mathbb{G}$}

To prove the main results of this paper, we first recall some necessary notions and remarks.
Firstly, we recall some preliminaries concerning stratified Lie groups (or so-called
Carnot groups). We refer the reader to  \cite{folland1982hardy,bonfiglioli2007stratified,stein1993harmonic}.

%------------------------
 \begin{definition}\label{def:stratified-Lie-algebra-yessir2019}
%----------------------------
  We say that  a Lie algebra  $\mathcal{G}$  is  stratified  if there is a direct sum vector space decomposition
%-----------------
\begin{align}\label{equ:lie-algebra-decomposition}
%-------
 \mathcal{G} =\oplus_{j=1}^{m} V_{j}  = V_{1} \oplus  \cdots \oplus  V_{m}
%----------
\end{align}
%------------
 such that  $\mathcal{G}$ is nilpotent of step $m$ if  $m$ is the smallest integer for which all Lie brackets (or iterated commutators) of order $m+1$ are zero, that is,
%-------------------------------
%-----------------
\begin{align*}
%-------
 [V_{1},V_{j}] =
%------------
 \begin{cases}
 %----------
 V_{j+1}, &  1\le j \le  m-1  \\
 0, & j\ge m
%------------
\end{cases}
%----------
\end{align*}
%------------
 holds.
%----------------------------
\end{definition}
%--------------------------

It is not difficult to find that the above $V_{1}$  generates the whole of the Lie algebra $\mathcal{G}$ by taking Lie brackets.

%--------------------------------
\begin{remark}\cite{zhu2003herz}   \label{rem:lie-algebra-decom-zhu2003herz}  %
%--------------------------------
Let $\mathcal{G} =  \mathcal{G}_{1}\supset  \mathcal{G}_{2} \supset \cdots \supset  \mathcal{G}_{m+1} =\{0\}$   denote the lower central series of  $\mathcal{G}$, and $\{X_{1},\dots,X_{N}\}$ be a basis for $V_{1}$ of $\mathcal{G}$.
%-------------------------------
\begin{enumerate}[leftmargin=2em,label=(\roman*),itemindent=1.5em]  %,itemindent=-0.3em
%--------------------------------------------------------------
\item  The direct sum decomposition  \labelcref{equ:lie-algebra-decomposition} can be constructed by identifying each $\mathcal{G}_{j}$ as a vector subspace of $\mathcal{G}$ and setting $ V_{m}=\mathcal{G}_{m}$ and $ V_{j}=\mathcal{G}_{j}\setminus \mathcal{G}_{j+1}$ for $j=1,\ldots,m-1$.
%-------------------------------------
\item   The dimension of  $\mathbb{G}$ at infinity as the integer $Q$ is given by
%-----------------
\begin{align*}
%-------
 Q = \sum_{j=1}^{m} j \dim(V_{j}) = \sum_{j=1}^{m} \dim(\mathcal{G}_{j}).
%----------
\end{align*}
%------------
\end{enumerate}
%-----------------------
\end{remark}
%------------------------

%------------------------
 \begin{definition}\label{def:stratified-Lie-group}
%----------------------------
     A Lie group $\mathbb{G}$ is said to be stratified when it is a connected simply-connected Lie group and its Lie algebra $\mathcal{G}$ is stratified.
%----------------------------
\end{definition}
%--------------------------

If $\mathbb{G}$ is stratified, then its Lie algebra $\mathcal{G}$ admits a canonical family of dilations $\{\delta_{r}\}$, namely, for $r>0$, $X_{k}\in V_{k}~(k=1,\ldots,m)$,
%-----------------
\begin{align*}
%-------
  \delta_{r} \Big( \sum_{k=1}^{m} X_{k} \Big)  =  \sum_{k=1}^{m} r^{k} X_{k},
%----------
\end{align*}
%------------
which are Lie algebra automorphisms.

By the Baker-Campbell-Hausdorff formula for sufficiently small elements $X$ and $Y$ of $\mathcal{G}$ one has
%------------------
\begin{align*}
%-------
 \exp X \exp Y=  \exp H(X,Y)= X+Y +\frac{1}{2}[X,Y]+\cdots
%----------
\end{align*}
%------------
where $\exp : \mathcal{G} \to \mathbb{G}$ is the exponential map, $H(X, Y )$ is an infinite linear
combination of $X$ and $Y$ and their Lie brackets, and the dots denote terms of order higher than two.

The following properties can be found in \cite{ruzhansky2019hardy}(see Proposition 1.1.1, or  Proposition 2.1 in \cite{yessirkegenov2019function} or Proposition 1.2 in \cite{folland1982hardy}).

%--------------------------------
\begin{proposition}\label{pro:2.1-yessirkegenov2019}
%------------------------------
 Let $\mathcal{G}$ be a nilpotent Lie algebra, and let $\mathbb{G}$ be the corresponding connected and simply-connected nilpotent Lie group. Then we have
%------------------------------
\begin{enumerate}[leftmargin=2em,label=(\arabic*),itemindent=1.5em]  %,itemindent=-0.3em
%--------------------------------
\item   The exponential map  $\exp: \mathcal{G} \to \mathbb{G}$  is a diffeomorphism. Furthermore, the group law $(x,y) \mapsto xy$ is a polynomial map if  $\mathbb{G}$ is identified with $\mathcal{G}$ via $\exp$.
%-------------------------------------
\item  If $\lambda$ is a Lebesgue measure on  $\mathcal{G}$, then $\exp\lambda$ is a bi-invariant Haar measure on  $\mathbb{G}$ (or a bi-invariant Haar measure $\mathd  x$ on  $\mathbb{G}$  is just the lift of Lebesgue measure on  $\mathcal{G}$ via $\exp$).
%-------------------------------
\end{enumerate}
%-----------------------------
\end{proposition}
%------------------------

\textbf{Notations:}
%-----------------------
\begin{itemize}
%-------------------------------------
 %  \item The number $Q= \sum_{j=1}^{m} j \dim(V_{j})$ denotes the homogenous dimension   of $\mathbb{G}$,
%-------------
   \item  $y^{-1}$ represents the inverse of $y\in \mathbb{G}$,
%-------------
   \item  $y^{-1}x$ stands for the group multiplication of $y^{-1}$  by $x$,
%-------------
  \item  Let the group identity element of $\mathbb{G}$  be referred to as the origin denotes by $e$,
%-----------
  \item $\chi_{E}$   denotes a characteristic function of a measurable set $E$ of $\mathbb{G}$,
%------------------
  \item $L^{p} ~(1\le p\le \infty)$  denotes the standard $L^{p} $-space with respect to the Haar measure $\mathd x$, with the $L^{p} $-norm $\|\cdot\|_{p}$.
%-------------
\end{itemize}
%--------------------

A homogenous norm on $\mathbb{G}$ is a continuous function $x\to \rho(x)$  from $\mathbb{G}$ to $[0,\infty)$, which is  $C^{\infty}$ on $\mathbb{G}\setminus\{0\}$ and satisfies
%-----------------
\begin{align*}
%-------
\begin{cases}
%------------
 \rho(x^{-1}) =  \rho(x), \\
 \rho(\delta_{t}x) =  t\rho(x) \ \ \text{for all}~  x \in \mathbb{G} ~\text{and}~ t > 0, \\
%------------
 \rho(e) =  0.
%------------
\end{cases}
%----------
\end{align*}
%------------
Moreover, there exists a constant $c_{0} \ge 1$ such that $\rho(xy) \le c_{0}(\rho(x) + \rho(y))$ for all $x,y \in \mathbb{G}$.

 With the norm above, we define the $\mathbb{G}$ ball centered at $x$ with radius $r$ by $B(x, r) = \{y \in \mathbb{G} : \rho(y^{-1}x) < r\}$,  and by $\lambda B$ denote the ball $B(x,\lambda r)$  with $\lambda>0$, let $B_{r} = B(e, r) = \{y \in \mathbb{G}  : \rho(y) < r\}$ be the open ball centered at $e$ with radius $r$,  which is the image under $\delta_{r}$ of $B(e, 1)$.
 And by $\sideset{^{\complement}}{}  {\mathop {B(x,r)}} = \mathbb{G}\setminus B(x,r)= \{y \in \mathbb{G} : \rho(y^{-1}x) \ge r\}$ denote the complement of $B(x, r)$.  Let  $\vert B(x,r)\vert $ be the Haar measure of the ball  $B(x,r)\subset \mathbb{G}$, and
 there exists $c_{1} =c_{1} (\mathbb{G})$ such that
%--------------------------------
\begin{align*}
%-------
 \vert B(x,r)\vert = c_{1} r^{Q}, \ \  \ \   x\in \mathbb{G}, r>0.
%----------
\end{align*}
%-----------
The most basic partial differential operator in a stratified Lie group is the sub-Laplacian associated with $X$ is the second-order partial differential operator on
 $\mathbb{G}$  given by
%--------------------------------
\begin{align*}
%-------
  \mathcal{L} =  \sum_{i=1}^{n} X_{i}^{2}.
%----------
\end{align*}
%-----------

%------------------------
\subsection{Maximal function}
%------------------------

Let $0\le \alpha <Q$ and $f: \mathbb{G} \to \mathbb{R}$ is a locally integrable function. The fractional
maximal function is defined by  %\cite{genebashvili1997weight}
%------------------
\begin{align*}
%-----------
   M_{\alpha}(f)(x)    &=  \sup_{B\ni x}  \dfrac{1}{\vert B\vert^{1-\alpha/Q}} \dint_{B}  \vert f(y)\vert \mathd y,
%-----------------
\end{align*}
%--------------
where the supremum is taken over all balls $B\subset \mathbb{G}$  containing $x$.

The fractional maximal function $ M_{\alpha}(f)$ coincides for $\alpha = 0$ with the Hardy-Littlewood maximal function $M(f)(x)\equiv M_{0}(f)(x)$.

The following propositions  can be found in \cite{kokilashvili1989fractional}.
%----------------
\begin{proposition} %(see \citet{kokilashvili1989fractional})
\label{pro:A-kokilashvili1989fractional}
%-------------------------------
Let $0\le\alpha<Q$ and $1< p < \gamma^{-1}=\frac{Q}{\alpha}$ with $\frac{1}{q}=\frac{1}{p}-\frac{\alpha}{Q}$. Then the following two conditions are equivalent:
%-------------------------------
\begin{enumerate}[leftmargin=2em,label=(\arabic*),itemindent=1.5em]
%--------------------------
\item  There is a constant $C>0$ such that for any  $f\in L_{\omega}^{p}(\mathbb{G})$  the inequality
%----------------
\begin{align*}
%-------
 \Big( \dint_{\mathbb{G}} \big( M_{\gamma}(f\omega^{\gamma}) (x) \big)^{q} \omega(x) \mathd x \Big)^{1/q} &\le C \Big(    \dint_{\mathbb{G}}\vert f(x) \vert^{p}  \omega(x) \mathd x \Big)^{1/p}
%----------
\end{align*}
%------------
holds.
%------------
\item       $\omega \in A_{1+q/p'}(\mathbb{G})$,   $p' = \frac{p}{p-1}$.
%------------
\end{enumerate}
%----------
\end{proposition}
%------------

%----------------
\begin{proposition} %(see \citet{kokilashvili1989fractional})
\label{pro:B-kokilashvili1989fractional}
%-------------------------------
Let   $0<\alpha<Q$, $\gamma=\alpha/Q$,  $q=(1-\gamma)^{-1}$, and $f\in L^{q}(\mathbb{G})$. Then the following two conditions are equivalent:
%-------------------------------
\begin{enumerate}[leftmargin=2em,label=(\arabic*),itemindent=1.5em]
%--------------------------
\item
%----------------
\begin{align*}
%-------
 \omega \{x\in \mathbb{G}: M_{\gamma}(f\omega^{\gamma})(x) >\lambda \}    &\le   C  \lambda^{-q} \Big(    \dint_{\mathbb{G}}\vert f(x) \vert   \mathd x \Big)^{q}
%----------
\end{align*}
%------------
with a constant $C>0$ independent of $f$ and $\lambda>0$.
%------------
\item       $\omega \in A_{1}(\mathbb{G})$.
%------------
\end{enumerate}
%----------
\end{proposition}
%------------

The following strong and weak-type boundedness of $M_{\alpha}$  can be obtained from \cref{pro:A-kokilashvili1989fractional,pro:B-kokilashvili1989fractional} when the weight $\omega=1$, see \citet{kokilashvili1989fractional} for more details.
And the first part  can also be obtained from \citet{bernardis1994two}.

%----------------
\begin{lemma}\label{lem:frac-maximal-kokilashvili1989fractional}
Let   $0<\alpha<Q$, $1\le p\le Q/\alpha$ with $1/q=1/p-\alpha/Q$, and $f\in L^{p}(\mathbb{G})$.
%-------------------------------
\begin{enumerate}[leftmargin=2em,label=(\arabic*),itemindent=1.5em]
%--------------------------
\item  If $1< p<Q/\alpha$,  then there exists a positive constant $C$ such that
%----------------
\begin{align*}
%-------
 \|M_{\alpha}(f)\|_{L^{q}(\mathbb{G})} &\le C  \|f\|_{L^{p}(\mathbb{G})}
%----------
\end{align*}
%------------
\item   If $p=1$,  then there exists a positive constant $C$ such that
%----------------
\begin{align*}
%-------
    \lvert \{x\in \mathbb{G}: M_{\alpha}(f)(x) >\lambda \} \rvert   &\le   C \big( \lambda^{-1}\|f\|_{L^{1}(\mathbb{G})}\big)^{Q/(Q-\alpha)}
%----------
\end{align*}
%------------
holds for all $\lambda>0$.
%------------
\end{enumerate}
%---------------------
%----------
\end{lemma}
%------------

\subsection{Lipschitz spaces on $\mathbb{G}$}

Next we give the definition of the Lipschitz spaces on $\mathbb{G}$, and state some basic properties and useful lemmas.

%--------------------------------
\begin{definition}[Lipschitz-type spaces on $\mathbb{G}$]   \label{def.lip-space} \
%--------------------------------排列环境--------------------------------
\begin{enumerate}[ label=(\arabic*),itemindent=1em]
%-------------------------------------
\item   Let $0<\beta <1$, we say a function $b$ belongs to the Lipschitz space $\dot{\Lambda}_{\beta}(\mathbb{G}) $ if there exists a constant $C>0$ such that for all  $x,y\in \mathbb{G}$,
%----------------
\begin{align} \label{inequ:lip-def-1}
%-------
\vert b(x)-b(y)\vert    &\le C(\rho(y^{-1}x))^{\beta},
%----------
\end{align}
%------------
where $\rho$ is the homogenous norm. The smallest such constant $C$ is called the $\dot{\Lambda}_{\beta}$  norm of $b$ and is denoted by $\|b\|_{\dot{\Lambda}_{\beta}(\mathbb{G})}$.
%-----------
    \label{enumerate:def-lip-1}
%-------------------------------------
\item (see \citet{macias1979lipschitz} ) Let $0<\beta <1$ and $1\le p<\infty$.   The space $\lip_{\beta,p}(\mathbb{G}) $ is defined to be the set of all locally integrable  functions $b$, i.e., there exists a positive constant $C $, such that
%----------------
\begin{align*}
%-------
      \sup_{B\ni x} \dfrac{1}{ \vert B \vert ^{\beta/Q}} \Big( \dfrac{1}{\lvert B\rvert }  \dint_{B} \lvert b(x)- b_{B}\rvert ^{p}\mathd x \Big)^{1/p} \le C,
%----------
\end{align*}
%------------
where the supremum is taken over every ball $B\subset \mathbb{G}$ containing $x$ and $b_{B}=\frac{1}{\vert B\vert } \int_{B} b(x) \mathd x$. The least constant $C$   satisfying the conditions above shall   be denoted by $\|b\|_{\lip_{\beta,p}(\mathbb{G})}$.
%-----------
    \label{enumerate:def-lip-2}
%-------------------------------------
\item (see \citet{macias1979lipschitz}) Let $0<\beta <1$. When  $  p=\infty$, we shall say that a locally integrable  functions $b$   belongs to  $\lip_{\beta,\infty}(\mathbb{G}) $  if there exists a constant $C$ such that
 %----------------
\begin{align*}
%-------
      \esssup_{x\in B} \dfrac{\lvert b(x)- b_{B}\rvert }{ \lvert B\rvert ^{\beta/Q}}   \le C
%----------
\end{align*}
%------------
holds  for every ball $B\subset \mathbb{G}$  with $b_{B}=\frac{1}{\vert B\vert } \int_{B} b(x) \mathd x$. And $\|b\|_{\lip_{\beta,\infty}(\mathbb{G})}$ stand for the least constant   $C$  satisfying the conditions above.
%-----------
    \label{enumerate:def-lip-3}
%-------------------
\end{enumerate}
%-------------------
\end{definition}
%--------------------

\begin{remark}  \label{rem.Lipschitz-def}
%-------------------------------
\begin{enumerate}[ label=(\roman*)]  %,itemindent=-0.3em,itemindent=1.5em
%--------------------------------------------------------------
\item  Similar to the definition of Lipschitz space $\dot{\Lambda}_{\beta}(\mathbb{G}) $ in \labelcref{enumerate:def-lip-1}, we also have the definition form as following  (see  \citet{krantz1982lipschitz,chen2010lipschitz,fan1995characterization} et al.)
%----------------
\begin{align*}
%-------
 \|b\|_{\dot{\Lambda}_{\beta}(\mathbb{G})}&= \sup_{x,y\in \mathbb{G}\atop y\neq e} \dfrac{\vert b(xy)- b(x)\vert }{(\rho(y))^{\beta}}   = \sup_{x,y\in \mathbb{G} \atop x\neq y} \dfrac{\vert b(x)-b(y)\vert }{(\rho(y^{-1}x))^{\beta}}.
%----------
\end{align*}
%------------
And $\|b\|_{\dot{\Lambda}_{\beta}(\mathbb{G})} =0$   if and only if $b$ is constant.
%-------------------------------------
\item  In \labelcref{enumerate:def-lip-2},  when   $p=1$, we have
%----------------
\begin{align*}
%-------
     \|b\|_{\lip_{\beta,1}(\mathbb{G})} =\sup_{B\ni x} \dfrac{1}{ \lvert B\rvert ^{\beta/Q}}\Big( \dfrac{1}{\lvert B\rvert }  \dint_{B} \lvert b(x)- b_{B}\rvert  \mathd x \Big) :=\|b\|_{\lip_{\beta}(\mathbb{G})}.
%----------
\end{align*}
%------------
\item    There are two basically different approaches to Lipschitz classes on the $n$-dimensional  Euclidean space. Lipschitz classes can be defined via Poisson (or Weierstrass) integrals  of $L^{p}$-functions, or, equivalently, by  means of higher order difference operators (see \citet{meda1988lipschitz}).
%------------
\end{enumerate}
%--------------------------
\end{remark}
%------------------------

%----------------
\begin{lemma} (see \cite{macias1979lipschitz,chen2010lipschitz,li2003lipschitz} ) \label{lem:2.2-li2003lipschitz}
Let   $0<\beta<1$ and the function $b(x)$ integrable on bounded subsets of $\mathbb{G}$.
%-------------------------------
\begin{enumerate}[leftmargin=2em,label=(\arabic*),itemindent=1.5em]  %,itemindent=-0.3em
%--------------------------
\item  When $1\le p<\infty$,  then
%----------------
\begin{align*}
%-------
 \|b\|_{\dot{\Lambda}_{\beta}(\mathbb{G})} &=  \|b\|_{\lip_{\beta}(\mathbb{G})} \approx  \|b\|_{\lip_{\beta,p}(\mathbb{G})}.
%----------
\end{align*}
%------------
\item   Let balls $B_{1}\subset B_{2}\subset \mathbb{G}$ and $b\in \lip_{\beta,p}(\mathbb{G})$ with $p\in [1,\infty]$. Then there exists a constant $C$ depends on $B_{1}$ and $B_{2}$ only, such that
%----------------
\begin{align*}
%-------
     \lvert b_{B_{1}}- b_{B_{2}} \rvert    &\le    C  \|b\|_{\lip_{\beta,p}(\mathbb{G})} \lvert B_{2}\rvert ^{\beta/Q}.
%----------
\end{align*}
%------------
\item   When $1\le p<\infty$, then there exists a constant $C$ depends on $\beta$ and $p$ only, such that
%----------------
\begin{align*}
%-------
    \lvert  b(x)-  b(y) \rvert    &\le   C  \|b\|_{\lip_{\beta,p}(\mathbb{G})} \lvert B\rvert ^{\beta/Q}
%----------
\end{align*}
%------------
holds for any ball $B$ containing $x$ and $y$.
%------------
\end{enumerate}
%---------------------
%----------
\end{lemma}
%------------

%------------------------
\subsection{Morrey spaces  on $\mathbb{G}$}
%------------------------

Morrey spaces were originally introduced by Morrey in \cite{morrey1938solutions} to study the local behavior of solutions to second-order elliptic partial differential equations.

%--------------------------------
\begin{definition}[Morrey-type spaces on $\mathbb{G}$\cite{eroglu2017characterizations}]   \label{def.morrey-space} \
%------------------------------
\begin{enumerate}[ label=(\arabic*),itemindent=1em]
%-------------------------------------
\item %(see \citet{eroglu2017characterizations})
Let $1\le p <\infty$ and $0\le \lambda \le Q$. The Morrey-type spaces $ L^{p,\lambda}(\mathbb{G})$ is defined by
%----------------
\begin{align*}
%-------
  L^{p,\lambda}(\mathbb{G})  &= \{ f\in L_{\loc}^{p}(\mathbb{G}): \|f\|_{L^{p,\lambda}(\mathbb{G})} < \infty \}
%----------
\end{align*}
%------------
with
%----------------
\begin{align*}
%-------
 \|f\|_{L^{p,\lambda}(\mathbb{G})}  &= \sup_{B\ni x \atop B\subset \mathbb{G}} \Big( \frac{1}{\lvert B\rvert ^{\lambda/Q}} \dint_{B} \vert f(y)\vert^{p} \mathd y \Big)^{1/p},
%----------
\end{align*}
%------------
where the supremum is taken over every ball $B\subset \mathbb{G}$ containing $x$.
%-----------
    \label{enumerate:def-morrey-1}
%-------------------------------------
\item  %(see \citet{eroglu2017characterizations})
Let $1\le p<\infty$ and $\varphi(x,r)$ be a positive measurable function on $\mathbb{G}\times (0,\infty)$.   The generalized Morrey space $ \mathcal{L}^{p,\varphi}(\mathbb{G})$ is defined for all functions $f\in L_{\loc}^{p}(\mathbb{G})$ by the finite norm
%----------------
\begin{align*}
%-------
     \|f\|_{\mathcal{L}^{p,\varphi}(\mathbb{G})}  &= \sup_{B\ni x \atop B\subset \mathbb{G}} \dfrac{1}{\varphi(x,r)} \Big( \frac{1}{\vert B\vert } \dint_{B} \vert f(y)\vert ^{p} \mathd y \Big)^{1/p},
%----------
\end{align*}
%------------
where the supremum is taken over every ball $B\subset \mathbb{G}$ containing $x$.
%-----------
    \label{enumerate:def-morrey-2}
%-------------------------------------
%----------
\end{enumerate}
%-------------------
\end{definition}
%--------------------

%--------------------
\begin{remark}[\citet{guliyev2020characterizations}]  \label{rem.morrrey-def}
%-------------------------------
\begin{enumerate}[ label=(\roman*)]  %,itemindent=-0.3em,itemindent=1.5em
%--------------------------------------------------------------
\item  It is well known that if   $1\le p <\infty$ then
%-----------------
\begin{align*}
%-------
  L^{p,\lambda}(\mathbb{G})  =
%-------
\begin{cases}
%------------
   L^{p}(\mathbb{G})  & \text{if}\ \lambda=0, \\
  L^{\infty}(\mathbb{G})  & \text{if}\ \lambda=Q,\\
  \Theta             & \text{if}\ \lambda<0 \ \text{or}\ \lambda>Q,
%------------
\end{cases}
%----------
\end{align*}
%------------
where $\Theta$ is the set of all functions equivalent to $0$ on $\mathbb{G}$.
%-------------------------------------
\item  In \labelcref{enumerate:def-morrey-2},  when   $1\le p<\infty $ %$1\le p\le q $, $\kappa \in (0,1)$
and  $0\le \lambda \le Q$, we have  $\mathcal{L}^{p,\varphi}(\mathbb{G})  = L^{p,\lambda}(\mathbb{G}) $
%-----------------
if $\varphi(x,r)=\vert B\vert ^{(\lambda/Q-1)/p}$ and $B\subset \mathbb{G}$ denotes the ball  with radius   $r$ and  containing $x$.
%------------
\end{enumerate}
%--------------------------
\end{remark}
%----------------------

We now recall the result on the boundedness of the fractional maximal operator  in the generalised Morrey spaces, which can be found in  \cite{guliyev2013boundedness} (theorem 3.2 and 3.3, see also \cite{nakai2006campanato}).

%----------------
\begin{proposition}[Spanne-type]  \label{pro:3.2-guliyev2013boundedness}
%---------------
 Let $1\le p<\infty$, $0\le\alpha <\frac{Q}{p}$, $\frac{1}{q} = \frac{1}{p} - \frac{\alpha}{Q}$ and $(\varphi_{1},\varphi_{2})$ satisfy the condition
%----------------
\begin{align*}
%-------
  \sup_{r<t<\infty} t^{\alpha-Q/p} \essinf_{t<s<\infty}  \varphi_{1}(x,s) s^{Q/p}   &\le C  \varphi_{2}(x,r),
%----------
\end{align*}
%------------
where  $C>0$ does not depend on $r$ and $x\in \mathbb{G}$.
%----------------
\begin{enumerate}[leftmargin=2em,label=(\arabic*),itemindent=1.5em]  %,itemindent=-0.3em
%--------------------------
\item  Then, for  $1< p <\infty$ and any $f\in \mathcal{L}^{p,\varphi_{1}}(\mathbb{G})$, there exists some positive  constant $C$ such that
%----------------
\begin{align*}
%-------
 \|M_{\alpha}f\|_{\mathcal{L}^{q,\varphi_{2}}(\mathbb{G})} &\le C  \|f\|_{\mathcal{L}^{p,\varphi_{1}}(\mathbb{G})}.
%----------
\end{align*}
%------------
\item Then, for  $p=1$ and any $f\in \mathcal{L}^{1,\varphi_{1}}(\mathbb{G})$, there exists some positive  constant $C$ such that
%----------------
\begin{align*}
%-------
     \|M_{\alpha}f\|_{W\mathcal{L}^{q,\varphi_{2}}(\mathbb{G})} &\le C  \|f\|_{\mathcal{L}^{1,\varphi_{1}}(\mathbb{G})}.
%----------
\end{align*}
%------------
\end{enumerate}
%----------
\end{proposition}
%------------

In the case $\alpha = 0$ and $p = q$, the the conclusions of \cref{pro:3.2-guliyev2013boundedness} are also valid.

%----------------
\begin{proposition}[Adams-typpe]  \label{pro:3.3-guliyev2013boundedness}
%---------------
 Let $1\le p< q<\infty$, $0<\alpha <\frac{Q}{p}$,  and  let $ \varphi(x,\tau)$ satisfy the condition
%----------------
\begin{align*}
%-------
  \sup_{r<t<\infty} t^{-Q} \essinf_{t<s<\infty}  \varphi(x,s) s^{Q}   &\le C  \varphi(x,r)
%------------
\\ \intertext{and}
%------------
     \sup_{r<t<\infty} t^{\alpha}   \varphi(x,\tau)^{1/p}   &\le C  r^{-\alpha p/(q-p)},
%----------
\end{align*}
%------------
where  $C>0$ does not depend on $r$ and $x\in \mathbb{G}$.
%----------------
\begin{enumerate}[leftmargin=2em,label=(\arabic*),itemindent=1.5em]  %,itemindent=-0.3em
%--------------------------
\item Then, for $1< p <\infty$ and any $f\in \mathcal{L}^{p,\varphi^{1/p}}(\mathbb{G})$, there exists some positive  constant $C$ such that
%----------------
\begin{align*}
%-------
 \|M_{\alpha}f\|_{\mathcal{L}^{q,\varphi^{1/q}}(\mathbb{G})} &\le C  \|f\|_{\mathcal{L}^{p,\varphi^{1/p}}(\mathbb{G})}.
%----------
\end{align*}
%------------
\item Then, for $p=1$ and any $f\in \mathcal{L}^{1,\varphi^{1/p}}(\mathbb{G})$, there exists some positive  constant $C$ such that
%----------------
\begin{align*}
%-------
     \|M_{\alpha}f\|_{W\mathcal{L}^{q,\varphi^{1/q}}(\mathbb{G})} &\le C  \|f\|_{\mathcal{L}^{1,\varphi}(\mathbb{G})}.
%----------
\end{align*}
%------------
\end{enumerate}
%----------
\end{proposition}
%------------

When $\varphi_{1}(x,r)=\vert B \vert ^{(\lambda/Q-1)/p}=c_{1}r^{Q(\lambda/Q-1)/p}$ and $\varphi_{2}(x,r)=\vert B \vert ^{(\mu/Q-1)/q}=c_{2}r^{Q(\mu/Q-1)/q}$,  we can summarize  the   results as follows   from   \cref{pro:3.2-guliyev2013boundedness,pro:3.3-guliyev2013boundedness}(see also Corollary 3.3 in \cite{guliyev2013boundedness}).

%----------------
\begin{lemma}    \label{lem:2.2-3.2-guliyev2013boundedness}
%---------------
 Let $0<\alpha <Q$, $1< p<Q/\alpha$  and $0<\lambda <Q-\alpha p$.
%----------------
\begin{enumerate}[leftmargin=2em,label=(\arabic*),itemindent=1.5em]  %,itemindent=-0.3em
%--------------------------
\item  If   $1/q=1/p-\alpha/(Q-\lambda)$, then there exists a positive  constant $C$ such that
%----------------
\begin{align*}
%-------
 \|M_{\alpha}f\|_{L^{q,\lambda}(\mathbb{G})} &\le C  \|f\|_{L^{p,\lambda}(\mathbb{G})}.
%----------
\end{align*}
%------------
for every $f\in L^{p,\lambda}(\mathbb{G})$.
%------------
\item  If   $1/q=1/p-\alpha/Q$ and $\lambda/p =\mu/q$.  Then there exists a positive  constant $C$ such that
%----------------
\begin{align*}
%-------
     \|M_{\alpha}f\|_{L^{q,\mu}(\mathbb{G})} &\le C  \|f\|_{L^{p,\lambda}(\mathbb{G})}
%----------
\end{align*}
%------------
for every $f\in L^{p,\lambda}(\mathbb{G})$.
%------------
\end{enumerate}
%----------
\end{lemma}
%------------

%------------------------
\section{Proofs of \cref{thm:lipschitz-main-1,thm:lipschitz-mb-main-2}} %Proof of the main results  the principal results
\label{sec:proof-mb}

We now give the proof of the main results. First, we prove \cref{thm:lipschitz-main-1}.

%--------------------

\begin{refproof}[Proof of \cref{thm:lipschitz-main-1}]
%--------------------------------
If  $b\in \dot{\Lambda}_{\beta}(\mathbb{G})$, then, using \labelcref{enumerate:def-lip-1} in \cref{def.lip-space}, we have
%-----------------
\begin{align}  \label{inequ:proof-mb-main-1-1}
%----------
\begin{aligned}
%-------
M_{b} (f)(x) &= \sup_{B\ni x} \vert B \vert ^{-1}  \dint_{B} \lvert b(x)-b(y)\rvert  \lvert f(y)\rvert  \mathd y  \\
&\le C\|b\|_{\dot{\Lambda}_{\beta}(\mathbb{G})} \sup_{B\ni x} \lvert B \rvert ^{-1}  \dint_{B} \lvert \rho(y^{-1}x)\rvert ^{\beta} \vert f(y) \vert  \mathd y     \\
&\le C\|b\|_{\dot{\Lambda}_{\beta}(\mathbb{G})} \sup_{B\ni x}  \dfrac{1}{\vert B \vert ^{1-\beta/Q}} \dint_{B}  \vert f(y) \vert  \mathd y        \\
&\le C\|b\|_{\dot{\Lambda}_{\beta}(\mathbb{G})}  M_{\beta} (f)(x).
%----------
\end{aligned}
%----------
\end{align}
%------------
Therefore, \labelcref{enumerate:thm-lip-main-1-2}, \labelcref{enumerate:thm-lip-main-1-3,enumerate:thm-lip-main-1-4} follow from  \cref{lem:frac-maximal-kokilashvili1989fractional} and above estimate.

\labelcref{enumerate:thm-lip-main-1-3} $\xLongrightarrow[]{\ \ \ \ }$ \labelcref{enumerate:thm-lip-main-1-1}:\ Suppose $ M_{b} $ is bounded from    $L^{p}(\mathbb{G})$ to $L^{q}(\mathbb{G})$ for some $p, q$ with $1<p<Q/\beta$ and $1/q = 1/p -\beta/Q$. For any ball $B\subset \mathbb{G}$ containing $x$, using the H\"{o}lder's inequality and noting that $1/p +1/q' =  1+\beta/Q$, one obtains
%-----------------
\begin{align*}
%-------
 \dfrac{1}{  \vert B \vert^{1+\beta/Q}}  \dint_{B}  \vert b(x)- b_{B} \vert \mathd x  &\le    \dfrac{1}{  \vert B \vert^{1+\beta/Q}}  \dint_{B} \Big(   \dfrac{1}{  \vert B \vert}  \dint_{B}  \vert b(x)- b(y) \vert \mathd y \Big)  \mathd x  \\
&=  \dfrac{1}{  \vert B \vert^{1+\beta/Q}}  \dint_{B} \Big(   \dfrac{1}{  \vert B \vert}  \dint_{B}  \vert b(x)- b(y) \vert\dchi_{B}(y) \mathd y \Big)  \mathd x  \\
&\le  \dfrac{1}{  \vert B \vert^{1+\beta/Q}}  \dint_{B}   M_{b}(\dchi_{B})(x)    \mathd x  \\
&\le  \dfrac{1}{  \vert B \vert^{1+\beta/Q}}  \Big(\dint_{B}  \big( M_{b}(\dchi_{B})(x) \big)^{q}   \mathd x  \Big)^{1/q}  \Big(\dint_{B}   \dchi_{B}(x)    \mathd x  \Big)^{1/q'}  \\
&\le  \dfrac{C}{  \vert B \vert^{1+\beta/Q}}   \|\dchi_{B}\|_{L^{p}(\mathbb{G})}   \| \dchi_{B}\|_{L^{q'}(\mathbb{G})}  \\
&\le C.
%----------
\end{align*}
%------------
This together with \cref{lem:2.2-li2003lipschitz} gives $b\in \dot{\Lambda}_{\beta}(\mathbb{G})$.

\labelcref{enumerate:thm-lip-main-1-4} $\xLongrightarrow[]{\ \ \ \ }$ \labelcref{enumerate:thm-lip-main-1-1}:\ Assume $ M_{b} $   satisfies the weak-type $(1,Q/(Q-\beta))$ estimates and \labelcref{inequ:lip-main-1-4} is true. In order to verify $b\in \dot{\Lambda}_{\beta}(\mathbb{G})$, for any fixed ball $B_{0}\subset \mathbb{G}$,  since for any $x\in B_{0}$,
%------------------
\begin{align*}
%-----------
   \vert b(x)-b_{B_{0}} \vert  &\le   \dfrac{1}{ \vert B_{0}\vert  }    \dint_{B_{0}}    \vert b(x)- b(y) \vert  \mathd y,
%-----------------
\end{align*}
%--------------
then, for all $x\in B_{0}$,
%------------------
\begin{align*}
%-----------
  M_{b}(\dchi_{B_{0}})(x)   &= \sup_{B\ni x}   \dfrac{1}{  \vert B \vert }    \dint_{B}    \vert b(x)- b(y) \vert\dchi_{B_{0}} (y) \mathd y  \\
    &\ge \dfrac{1}{ \vert B_{0}\vert  }    \dint_{B_{0}}    \vert b(x)- b(y) \vert\dchi_{B_{0}} (y) \mathd y   \\
    &= \dfrac{1}{ \vert B_{0}\vert  }    \dint_{B_{0}}    \vert b(x)- b(y) \vert  \mathd y   \\
    &\ge  \vert b(x)-b_{B_{0}} \vert.
%-----------------
\end{align*}
%--------------

Thus, combined with \labelcref{inequ:lip-main-1-4}, we have
%------------------
\begin{align*}
%-----------
  \big\vert \{x\in B_{0}:  \vert b(x)-b_{B_{0}} \vert >\lambda\} \big\vert  &\le \big\vert \{x\in B_{0}: M_{b}(\dchi_{B_{0}})(x)  >\lambda\} \big\vert    \\
  &\le C \left( \lambda^{-1}    \|\dchi_{B_{0}}  \|_{L^{1}(\mathbb{G})} \right)^{Q/(Q-\beta)}  \\
    &\le C \left( \lambda^{-1}    \vert  B_{0} \vert   \right)^{Q/(Q-\beta)}.
%-----------------
\end{align*}
%--------------

Set $t>0$  be a constant to be determined later, then applying Fubini's theorem, one get
%------------------
\begin{align*}
%-----------
  \dint_{B_{0}}    \vert b(x)-b_{B_{0}} \vert  \mathd x  &= \dint_{0}^{\infty}  \big\vert \{x\in B_{0}:  \vert b(x)-b_{B_{0}} \vert >\lambda\} \big\vert  \mathd\lambda \\
   &= \dint_{0}^{t}  \big\vert \{x\in B_{0}:  \vert b(x)-b_{B_{0}} \vert >\lambda\} \big\vert  \mathd\lambda \\
   &\qquad +   \dint_{t}^{\infty}  \big\vert \{x\in B_{0}:  \vert b(x)-b_{B_{0}} \vert >\lambda\} \big\vert  \mathd\lambda \\
  &\le t\vert B_{0}\vert  + C \dint_{t}^{\infty}  \big( \lambda^{-1}    \vert  B_{0} \vert   \big)^{Q/(Q-\beta)} \mathd\lambda \\
  &\le t\vert B_{0}\vert  + C \vert  B_{0} \vert ^{Q/(Q-\beta)}  \dint_{t}^{\infty}    \lambda^{-Q/(Q-\beta)} \mathd\lambda \\
    &\le C \left(  t\vert B_{0}\vert +    \vert  B_{0} \vert ^{Q/(Q-\beta)} t^{1-Q/(Q-\beta)} \right).
%-----------------
\end{align*}
%--------------
Let $t=\vert  B_{0} \vert ^{\beta/Q}$  in the above estimate, we get
%------------------
\begin{align*}
%-----------
  \dint_{B_{0}}    \vert b(x)-b_{B_{0}} \vert  \mathd x  &\le C   \vert B_{0} \vert^{1+ \beta/Q}.
%-----------------
\end{align*}
%--------------
It follows from \cref{lem:2.2-li2003lipschitz}  that $b\in \dot{\Lambda}_{\beta}(\mathbb{G})$ since $B_{0}$ is an arbitrary ball in $\mathbb{G}$.

The proof of \cref{thm:lipschitz-main-1} is completed since \labelcref{enumerate:thm-lip-main-1-2} $\xLongrightarrow[]{\ \ \ \ }$ \labelcref{enumerate:thm-lip-main-1-1} follows from \labelcref{enumerate:thm-lip-main-1-3} $\xLongrightarrow[]{\ \ \ \ }$ \labelcref{enumerate:thm-lip-main-1-1}.
%--------------------------------
\end{refproof}
%-----------------------

%--------------------
\begin{refproof}[Proof of \cref{thm:lipschitz-mb-main-2}]
%----------------------
\labelcref{enumerate:thm-lip-mb-main-2-1} :\   We first prove that  the necessary condition. Assume $b\in \dot{\Lambda}_{\beta}(\mathbb{G})$, using \labelcref{inequ:proof-mb-main-1-1} and   \cref{lem:2.2-3.2-guliyev2013boundedness}, we obtain
%-----------------
\begin{align*}
%-------
\|M_{b} (f)\|_{L^{q,\lambda}(\mathbb{G})}  &\le C\|b\|_{\dot{\Lambda}_{\beta}(\mathbb{G})} \|M_{\beta}f\|_{L^{q,\lambda}(\mathbb{G})}  \le C \|b\|_{\dot{\Lambda}_{\beta}(\mathbb{G})}  \|f\|_{L^{p,\lambda}(\mathbb{G})}.
%----------
\end{align*}
%------------

  We now prove that the sufficient condition. If $ M_{b} $  is bounded from  $L^{p,\lambda}(\mathbb{G})$ to $L^{q,\lambda}(\mathbb{G})$, then for any ball $B\subset \mathbb{G}$,
%-----------------
\begin{align*}
%-------
\vert B \vert^{-\beta/Q}  \Big(   \vert B \vert^{-1}  \dint_{B}  \vert b(x)- b_{B} \vert^{q}  \mathd x \Big)^{1/q}
% &\le  \vert B \vert^{-\beta/Q}  \bigg(   \vert B \vert^{-1}  \dint_{B} \Big(   \vert B \vert^{-1}  \dint_{B}  \vert b(x)- b(y) \vert \dchi_{B}(y) \mathd y \Big)^{q} \mathd x \bigg)^{1/q}   \\
 &\le  \vert B \vert^{-\beta/Q}  \Big(   \vert B \vert^{-1}  \dint_{B} \big( M_{b}(\dchi_{B})(x)\big)^{q}  \mathd x \Big)^{1/q}  \\
%&=  \vert B \vert^{-\beta/Q}  \vert B \vert^{(\lambda/Q-1)/q} \Big(   \vert B \vert^{-\lambda/Q}  \dint_{B} \big( M_{b}(\dchi_{B})(x)\big)^{q}  \mathd x \Big)^{1/q}  \\
 &\le    \vert B \vert^{-\beta/Q-1/q+\lambda/(Qq)}    \|M_{b}(\dchi_{B})\|_{L^{q,\lambda}(\mathbb{G})}    \\
 &\le C   \vert B \vert^{-\beta/Q-1/q+\lambda/(Qq)}     \|\dchi_{B}\|_{L^{p,\lambda}(\mathbb{G})}   \\
  &\le  C,
%----------
\end{align*}
%------------
where in the last step we have used $1/q=1/p-\beta/(Q-\lambda)$ and the fact
%-----------------
\begin{align}   \label{inequ:proof-mb-main-2-3}
%-------
  \|\dchi_{B}\|_{L^{p,\lambda}(\mathbb{G})} &\le  \vert B \vert^{(1-\lambda/Q)/p}.
%----------
\end{align}
%------------

It follows from \cref{lem:2.2-li2003lipschitz} that $b\in \dot{\Lambda}_{\beta}(\mathbb{G})$. This completes the proof.

\labelcref{enumerate:thm-lip-mb-main-2-2} :\  By a similar proof to \labelcref{enumerate:thm-lip-mb-main-2-1} in \cref{thm:lipschitz-mb-main-2}, we can obtain the desired result.
%-----------------

%--------------------
\end{refproof}
%-----------------------

%------------------------
\section{Proofs of \cref{thm:lipschitz-nonlinear-main-1,thm:lipschitz-nonlinear-main-3,thm:lipschitz-nonlinear-main-4}}
 \label{sec:proof-nonlinear}

%-----------------------

 Now, we prove \cref{thm:lipschitz-nonlinear-main-1}.

%--------------------
\begin{refproof}[Proof of \cref{thm:lipschitz-nonlinear-main-1}]
%----------------------
\labelcref{enumerate:thm-lip-nonlinear-main-1-1} $\xLongrightarrow[]{\ \ \ \ }$ \labelcref{enumerate:thm-lip-nonlinear-main-1-2}:\  For any fixed $x \in \mathbb{G}$ such that $M(f)(x) <\infty$, since $b \geq 0$ then
%-----------------
\begin{align}  \label{inequ:proof-nonlinear-main-2-1}
%-------
\begin{aligned}
%-------
\vert [b,M] (f)(x)\vert  &= \vert b(x)M(f)(x)-M(bf)(x)\vert   \\
% &= \bigg\vert  \sup_{B\ni x}  \vert B \vert^{-1}  \dint_{B} b(x) \vert f(y)  \vert  \mathd y -\sup_{B\ni x}  \vert B \vert^{-1}  \dint_{B} b(y) \vert f(y)  \vert  \mathd y  \bigg\vert  \\
 &\le \sup_{B\ni x}  \vert B \vert^{-1}  \dint_{B}  \vert b(x)- b(y) \vert  \vert f(y)  \vert \mathd y  \\
&=  M_{b} (f)(x).
%----------
 \end{aligned}
%-------
\end{align}
%------------
It follows from  \cref{thm:lipschitz-main-1} that $[b, M]$ is bounded from  $L^{p}(\mathbb{G})$ to $L^{q}(\mathbb{G})$ since $b\in \dot{\Lambda}_{\beta}(\mathbb{G})$.

\labelcref{enumerate:thm-lip-nonlinear-main-1-2} $\xLongrightarrow[]{\ \ \ \ }$ \labelcref{enumerate:thm-lip-nonlinear-main-1-3}:\  For any fixed ball $B\subset \mathbb{G}$ and all $x \in B$, one have
%-----------------
\begin{align*}
%-------
M (\dchi_{B})(x) =  \dchi_{B}(x)  \ \ \text{and} \ \  M (b\dchi_{B})(x)  =M_{B} (b)(x).
%----------
\end{align*}
%------------
Then
%-----------------
\begin{align}\label{inequ:proof-nonlinear-main-2-2}
%-------
\begin{aligned}
%-------
&  \vert B \vert^{-\beta/Q}  \Big(   \vert B \vert^{-1}  \dint_{B}  \vert b(x)-M_{B} (b)(x) \vert^{q}  \mathd x \Big)^{1/q}  \\
%&=  \vert B \vert^{-\beta/Q}  \Big(   \vert B \vert^{-1}  \dint_{B} \vert  b(x)M (\dchi_{B})(x) -M (b\dchi_{B})(x)\vert ^{q}  \mathd x \Big)^{1/q}   \\
 &=  \vert B \vert^{-\beta/Q}  \Big(   \vert B \vert^{-1}  \dint_{B}  \vert  [b,M] (\dchi_{B})(x) \vert^{q}  \mathd x \Big)^{1/q}  \\
 &\le    \vert B \vert^{-\beta/Q-1/q} \|[b,M] (\dchi_{B})\|_{L^{q}(\mathbb{G})}    \\
  &\le  C  \vert B \vert^{-\beta/Q-1/q} \| \dchi_{B} \|_{L^{p}(\mathbb{G})}    \le C,
%----------
 \end{aligned}
%----------
\end{align}
%------------
which implies  \labelcref{enumerate:thm-lip-nonlinear-main-1-3} since the ball $B\subset \mathbb{G}$ is arbitrary.

\labelcref{enumerate:thm-lip-nonlinear-main-1-3} $\xLongrightarrow[]{\ \ \ \ }$ \labelcref{enumerate:thm-lip-nonlinear-main-1-1}:\
To prove $b\in \dot{\Lambda}_{\beta}(\mathbb{G})$, by  \cref{lem:2.2-li2003lipschitz}, it suffices to verify that there is a constant $C>0$ such that for all balls $B\subset \mathbb{G}$, one get
%-----------------
\begin{align} \label{inequ:proof-nonlinear-main-1-31-1}
%-------
 \vert B \vert^{-1-\beta/Q}   \dint_{B} \vert  b(x)-b_{B} \vert   \mathd x   \le C.
%----------
\end{align}
%------------

For any fixed ball $B\subset \mathbb{G}$, let $E=\{x \in B:b(x) \le b_{B}\}$ and $F=\{x \in B:b(x) > b_{B}\}$. The following equality is trivially true (modifying the argument in \cite{bastero2000commutators}, page 3331):
%-----------------
\begin{align*}
%-------
   \dint_{E}  \vert b(x)- b_{B} \vert  \mathd x   =    \dint_{F}  \vert b(x)- b_{B} \vert  \mathd x.
%----------
\end{align*}
%------------
Since for any $x\in E$ we have $b(x) \le b_{B} \le M_{B}(b)(x)$, then for any $x\in E$,
%-----------------
\begin{align*}
%-------
     \vert b(x)- b_{B} \vert   \le      \vert b(x)-M_{B} (b)(x) \vert.
%----------
\end{align*}
%------------
Therefore
%-----------------
\begin{align} \label{inequ:proof-nonlinear-main-1-31-2}
\begin{aligned}
%-------
  \dfrac{1}{ \vert B \vert^{1+\beta/Q} } \dint_{B}  \vert b(x)- b_{B} \vert  \mathd x  &=  \dfrac{1}{ \vert B \vert^{1+\beta/Q} } \dint_{E\cup F}  \vert b(x)- b_{B} \vert  \mathd x  \\
  &=  \dfrac{2}{ \vert B \vert^{1+\beta/Q} } \dint_{E}  \vert b(x)- b_{B} \vert  \mathd x  \\
  &\le  \dfrac{2}{ \vert B \vert^{1+\beta/Q} } \dint_{E}  \vert b(x)-M_{B} (b)(x) \vert  \mathd x  \\
    &\le  \dfrac{2}{ \vert B \vert^{1+\beta/Q} } \dint_{B}  \vert b(x)-M_{B} (b)(x) \vert  \mathd x.
%----------
\end{aligned}
\end{align}
%------------

On the other hand, it follows from H\"{o}lder's inequality and \labelcref{inequ:lip-nonlinear-main-1-3} that
%-----------------
\begin{align*}
%-------
 &\ \dfrac{1}{ \vert B \vert^{1+\beta/Q} } \dint_{B}   \vert b(x)-M_{B} (b)(x) \vert  \mathd x  \\
  &\le  \dfrac{1}{ \vert B \vert^{1+\beta/Q} } \bigg( \dint_{B}   \vert b(x)-M_{B} (b)(x) \vert^{q}  \mathd x \bigg)^{1/q}   \vert B \vert^{1/q'}\\
  &\le  \dfrac{1}{ \vert B \vert^{\beta/Q} } \bigg(  \vert B \vert^{-1} \dint_{B}   \vert b(x)-M_{B} (b)(x) \vert^{q}  \mathd x \bigg)^{1/q} \\
    &\le  C.
%----------
\end{align*}
%------------
This together with \labelcref{inequ:proof-nonlinear-main-1-31-2} gives \labelcref{inequ:proof-nonlinear-main-1-31-1}, and so we achieve $b\in \dot{\Lambda}_{\beta}(\mathbb{G})$.

In order to prove $b \ge 0$, it suffices to show $b^{-}=0$, where $b^{-}=- \min\{b, 0\}$. Let $b^{+}= \vert b  \vert-b^{-}$, then $b=b^{+}-b^{-}$. For any fixed ball $B\subset \mathbb{G}$, observe that
%------------------
\begin{align*}
%-----------
 0\le b^{+}(x) \le  \vert b(x)  \vert \le   M_{B}(b)(x)
%-----------------
\end{align*}
%-----------
for $x\in B$ and thus we have that, for $x\in B$,
%------------------
\begin{align*}
%-----------
 0\le b^{-}(x) \le  M_{B}(b)(x) -b^{+}(x) \le  M_{B}(b)(x) -b^{+}(x)+b^{-}(x) =M_{B}(b)(x) -b(x).
%-----------------
\end{align*}
%-----------
Then, it follows from \labelcref{inequ:lip-nonlinear-main-1-3} that, for any ball $B\subset \mathbb{G}$,
%-----------------
\begin{align*}
%-------
  \dfrac{1}{ \vert B \vert} \dint_{B}  b^{-}(x)  \mathd x  &\le    \dfrac{1}{ \vert B \vert} \dint_{B}  \vert M_{B}(b)(x) -b(x)\vert   \mathd x   \\
  &\le    \bigg( \dfrac{1}{ \vert B \vert} \dint_{B}   \vert b(x)-M_{B} (b)(x) \vert^{q}  \mathd x \bigg)^{1/q} \\
  &=   \vert B \vert^{\beta/Q}  \bigg( \dfrac{1}{ \vert B \vert^{\beta/Q} } \Big( \dfrac{1}{ \vert B \vert} \dint_{B}   \vert b(x)-M_{B} (b)(x) \vert^{q}  \mathd x \Big)^{1/q} \bigg) \\
  &\le C   \vert B \vert^{\beta/Q}.
%----------
\end{align*}
%------------
Thus, $b^{-}=0$ follows from Lebesgue's differentiation theorem.

The proof of \cref{thm:lipschitz-nonlinear-main-1} is completed.

%--------------------
\end{refproof}
%-----------------------

%--------------------
\begin{refproof}[Proof of \cref{thm:lipschitz-nonlinear-main-3}]
%----------------------
\labelcref{enumerate:thm-lip-nonlinear-main-3-1} $\xLongrightarrow[]{\ \  }$ \labelcref{enumerate:thm-lip-nonlinear-main-3-2}:\   We first prove that  the necessary condition. Assume $b\in \dot{\Lambda}_{\beta}(\mathbb{G})$  and $b\ge 0$. Using \labelcref{inequ:proof-nonlinear-main-2-1} and  \labelcref{enumerate:thm-lip-mb-main-2-1} in \cref{thm:lipschitz-mb-main-2}, it  is not difficult to find that $[b,M ]$ is bounded from  $L^{p,\lambda}(\mathbb{G})$ to $L^{q,\lambda}(\mathbb{G})$.
%-----------------

%----------
 \labelcref{enumerate:thm-lip-nonlinear-main-3-2} $\xLongrightarrow[]{\ \  }$  \labelcref{enumerate:thm-lip-nonlinear-main-3-1}:\
%------------
  We now prove that the sufficient condition. Assume that $[b,M ]$ is bounded from  $L^{p,\lambda}(\mathbb{G})$ to $L^{q,\lambda}(\mathbb{G})$. Similarly to \labelcref{inequ:proof-nonlinear-main-2-2}, for any  ball $B\subset \mathbb{G}$, we obtain
%-----------------
\begin{align*}
%-------
\begin{aligned}
%-------
  &\  \vert B \vert^{-\beta/Q}  \Big(   \vert B \vert^{-1}  \dint_{B}  \vert b(x)-M_{B} (b)(x) \vert^{q}  \mathd x \Big)^{1/q}    \\
 &=  \vert B \vert^{-\beta/Q}  \Big(   \vert B \vert^{-1}  \dint_{B}  \vert  [b,M] (\dchi_{B})(x) \vert^{q}  \mathd x \Big)^{1/q}  \\
 &\le    \vert B \vert^{\lambda/(Qq)-\beta/Q-1/q} \|[b,M] (\dchi_{B})\|_{L^{q,\lambda}(\mathbb{G})}    \\
  &\le  C  \vert B \vert^{\lambda/(Qq)-\beta/Q-1/q} \| \dchi_{B} \|_{L^{p,\lambda}(\mathbb{G})}    \le C,
%----------
 \end{aligned}
%----------
\end{align*}
%------------
%------------
where in the last step we have used $1/q=1/p-\beta/(Q-\lambda)$ and \labelcref{inequ:proof-mb-main-2-3}.
%-----------------

 Using \cref{thm:lipschitz-nonlinear-main-1}, we can obtain  that $b\in \dot{\Lambda}_{\beta}(\mathbb{G})$ and $b\ge 0$.

%--------------------
\end{refproof}
%-----------------------

%-----------------------

%--------------------
\begin{refproof}[Proof of \cref{thm:lipschitz-nonlinear-main-4}]
%----------------------
By the same way of the proof of \cref{thm:lipschitz-nonlinear-main-3}, \cref{thm:lipschitz-nonlinear-main-4} can be proven. We omit the details.

%--------------------
\end{refproof}
%-----------------------

\section*{Declarations}

\bmhead{Acknowledgments}

%The authors are greatly indebted to the referees for useful comments.
This work is supported partly  by  the National Natural Science Foundation of China (Grant No.11571160), Scientific Project-HLJ (No.2019-
KYYWF-0909) and Youth Project-HLJ (No.2020YQ07).

\bmhead{Competing interests}

The authors declare that they have no competing interests.

 \bmhead{Data Availability Statement}
My manuscript has no associate data.

 \bmhead{Authors' contributions}
%All authors contributed equally and significantly in writing this paper. All authors read and approved the final manuscript.

The first draft of the manuscript was written by Jianglong Wu and all authors commented on previous versions of the manuscript. All authors read and approved the final manuscript

 \bibliographystyle{bst/sn-basic}  %  (不好 sn-chicago 顺序颠倒) 姓前  sn-basic 按字母顺序，chicago elsarticle-harv 可以用\citet, 年份太靠前
%\bibliographystyle{bst/sn-aps}  % 姓后  按引用顺序  无法用 \citet
%\bibliographystyle{elsarticle-num-names}   % (For authoryear Elsevier citations), 可以 用\citet
%\bibliographystyle{unsrtnat}              % plainnat 全称，姓后，按字母 unsrtnat 按引用
% \bibliographystyle{chicago}   % apalike  acm 姓前 全大写 缩写名  按字母排序  无法用 \citet
%-------------- %

\bibliography{wu-reference}


\end{document}